\documentclass[a4paper]{amsart}
\usepackage[all]{xy}
\usepackage{hyperref}
\usepackage{enumerate}

\usepackage{amsmath}
\usepackage{amsthm}
\usepackage{amssymb}
\usepackage{amscd}
\usepackage{graphicx}
\usepackage{epsfig}
\usepackage[all]{xy}

\numberwithin{equation}{section}
\allowdisplaybreaks

\newtheorem{theorem}{Theorem}[section]
\newtheorem{lemma}{Lemma}[section]
\newtheorem{proposition}{Proposition}[section]
\newtheorem{corollary}{Corollary}[section]

\theoremstyle{definition}
\newtheorem{definition}{Definition}[section]
\newtheorem{example}{Example} [section]

\newtheorem{remark}{Remark}[section]

\begin{document}
\title{A Note on Proper Poisson Actions}
\author{Rui Loja Fernandes}
\address{Depart.~de Matem\'{a}tica, 
Instituto Superior T\'{e}cnico, 1049-001 Lisboa, PORTUGAL} 
\email{rfern@math.ist.utl.pt}
\thanks{Supported in part by FCT/POCTI/FEDER and by grant
  POCTI/MAT/57888/2004.}

\begin{abstract}
We show that the fixed point set of a proper action of a Lie
group $G$ on a Poisson manifold $M$ by Poisson automorphisms has a
natural induced Poisson structure and we give several applications.
\end{abstract}
\maketitle

\section{Introduction}

{In} the present work, we consider a \emph{Poisson action} $G\times M\to M$ 
of a Lie group $G$ on a Poisson manifold $M$: this means that each
element $g\in G$ acts by a Poisson diffeomorphism of $M$. We recall
that the action is called \emph{proper} if the map:
\[ G\times M\to M\times M,\quad (g,p)\mapsto (p,g\cdot p),\]
is a proper map\footnote{A map $f:X\to Y$ between two topological
spaces is called \emph{proper} if for every compact subset $K\subset Y$, 
the inverse image $f^{-1}(K)$ is compact.}.  As usual, we will denote by
$M^G$ the \emph{fixed point set} of the action:
\[ M^G=\{p\in M: g\cdot p=p, \forall g\in G\}.\]
For proper actions, the connected components of the fixed point set
$M^G$ are (embedded) submanifolds of $M$. Notice that these components
may have different dimensions.

The main result of this paper is the following: 

\begin{theorem}
\label{Fernandes:MainTheorem}
Let $G\times M\to M$ be a proper Poisson action. Then the fixed point
set $M^G$ has a natural induced Poisson structure.
\end{theorem}

This result is a generalization to Poisson geometry of a well-known
proposition in symplectic geometry, due to Guillemin and Sternberg (see
\cite{Fernandes:GuSt}, Theorem 3.5), stating that fixed point sets of
symplectic actions are symplectic submanifolds. We stress that the
fixed point set \emph{is not} a Poisson submanifold. This happens already
in the symplectic case. In the general Poisson case, $M^G$ will be a
Poisson-Dirac submanifold in the sense of Crainic and Fernandes
(see \cite{Fernandes:CrFe}, Section 8) and Xu (\cite{Fernandes:Xu}).

Proper symplectic/Poisson actions have been study intensively in the
last 15 years. For example, the theory of (singular) reduction for
Hamiltonian systems has been developed extensively for these kind of
actions. We refer the reader to the recent monograph by Ortega and
Ratiu \cite{Fernandes:OrRa} for a nice survey of results in this
area. Theorem \ref{Fernandes:MainTheorem} should have important
applications in symmetry reduction, and this is one of our main
motivations for this work. We refer the reader for an upcoming
publication (\cite{Fernandes:FeOrRa}).

This paper is organized as follows. In Section 1, we recall the notion
of a Poisson-Dirac submanifold, and some related results which are needed for
the proof of Theorem \ref{Fernandes:MainTheorem}. In Section 2, we prove our
main result. In Section 3, we deduce some consequences and give some
applications.

\section{Poisson-Dirac submanifolds}

Let $M$ be a Poisson manifold. For background in Poisson geometry we
refer the reader to Vaisman's book \cite{FernandesVai}. We will
denote by $\pi\in{\mathfrak X}^2(M)$ the Poisson bivector field so that the
Poisson bracket is given by:
\[ \{f,g\}=\pi({\mathrm d}f,{\mathrm d}g),\qquad \forall f,g\in C^\infty(M).\] 
Recall that a Poisson submanifold $N\subset M$ is a submanifold which
has a Poisson bracket and for which the inclusion
$i:N\hookrightarrow M$ is a Poisson map:
\[ 
\{f\circ i,g\circ i\}_M=\{f,g\}_N\circ i,\qquad \forall f,g\in C^\infty(N).
\]
Such Poisson submanifolds are, in a sense, extremely rare. In fact,
they are collections of open subsets of symplectic leaves of $M$. 

\begin{example}
Let $M$ be a symplectic manifold with symplectic form $\omega$. Recall
that a \emph{symplectic submanifold} is a submanifold
$i:N\hookrightarrow M$ such that the restriction $i^*\omega$ is a
symplectic form on $N$. For every even dimension $0\le 2i\le \dim M$
there are symplectic submanifolds of dimension $2i$. On the other hand,
the only Poisson submanifolds are the open subsets of $M$. 
\end{example}

Crainic and Fernandes in \cite{Fernandes:CrFe} introduce the following
natural extension of the notion of a Poisson submanifold:

\begin{definition}
\label{Fernandes:defnDiracSub}
Let $M$ be a Poisson manifold. A submanifold $N\subset M$ is called a
\textbf{Poisson-Dirac submanifold} if $N$ is a Poisson manifold such that: 
\begin{enumerate}
\item[(i)] the symplectic foliation of $N$ is $N\cap{\mathcal F}=\{L\cap
    N:L\in{\mathcal F}\}$, and 
\item[(ii)] for every leaf $L\in{\mathcal F}$, $L\cap N$ is a symplectic
  submanifold of $L$.
\end{enumerate}
\end{definition}

Note that if $(M,\{\cdot,\cdot\})$ is a Poisson manifold, then the
symplectic foliation with the induced symplectic forms on the leaves,
gives a smooth (singular) foliation with a smooth family of symplectic
forms. Conversely, given a manifold $M$ with a foliation ${\mathcal F}$
furnished with a smooth family of symplectic forms on the leaves,
then we have a Poisson bracket on $M$ defined by the
formula\footnote{In a Poisson (or symplectic) manifold, we will denote
by $X_f$ the Hamiltonian vector field associated with a function
$f:M\to{\mathbb R}$.}
\[ \{f,g\}\equiv X_f(g),\]
for which the associated symplectic foliation is precisely
${\mathcal F}$. Hence, a Poisson structure can be defined by specifying its
symplectic foliation.  It follows that a submanifold $N$ of a Poisson
manifold $M$ has at most one Poisson structure satisfying conditions
(i) and (ii) above, and this Poisson structure is completely
determined by the Poisson structure of $M$.

\begin{example}
If $M$ is a symplectic manifold, then there is only one symplectic
leave, and the Poisson-Dirac submanifolds are precisely the symplectic
submanifolds of $M$.
\end{example}

Therefore, we see that the notion of a Poisson-Dirac submanifold
generalizes to the Poisson category the notion of a symplectic submanifold.

\begin{example}
Let $L$ be a symplectic leaf of a Poisson manifold, and $N\subset M$ a
submanifold which is transverse to $L$ at some $x_0$:
\[ T_{x_0}M=T_{x_0}L\oplus T_{x_0}N.\]
Then one can check that conditions (i) and (ii) in Definition
\ref{Fernandes:defnDiracSub} are satisfied in some open subset in $N$
containing $x_0$. In other words, if $N$ is small enough then it is a
Poisson-Dirac submanifold. Sometimes one calls the Poisson structure
on $N$ the transverse Poisson structure to $L$ at $x_0$ (up to Poisson
diffeomorphisms, this structure does not depend on the transversal
$N$).
\end{example}

The two conditions in Definition \ref{Fernandes:defnDiracSub} are not
very practical to use. Let us give some alternative criteria to
determine if a given submanifold is a Poisson-Dirac submanifold.

Observe that condition (ii) in the definition means that the
symplectic forms on a leaf $L\cap N$ are the pull-backs $i^*\omega_L$,
where $i:N\cap L\hookrightarrow L$ is the inclusion into a leaf and
$\omega_L\in\Omega^2(L)$ is the symplectic form. Denoting by
$\#:T^*M\to TM$ the bundle map determined by the Poisson bivector
field, we conclude that we must have\footnote{For a subspace
  $W$ of a vector space $V$, we denote by $W^0\subset V^*$ its
  annihilator. Similarly, for a vector subbundle $E\subset F$, we 
  denote by $E^0\subset F^*$ its annihilator subbundle.}:
\begin{equation}
\label{Fernandes:eqDirac}
TN\cap \#(TN^0)=\{0\},
\end{equation}
since the left-hand side is the kernel of the pull-back
$i^*\omega_L$. If this condition holds, then at each point $x\in N$ we
obtain a bivector $\pi_N(x)\in\wedge^2 T_xN$, and one can prove (see
\cite{Fernandes:CrFe}): 

\begin{proposition}
\label{Fernandes:propDirac}
Let $N$ be a submanifold of a Poisson manifold $M$, such that
\begin{enumerate}
\item[(a)] equation (\ref{Fernandes:eqDirac}) holds, and
\item[(b)] the induced tensor $\pi_N$ is smooth.
\end{enumerate}
Then $\pi_N$ is a Poisson tensor and $N$ is a Poisson-Dirac submanifold.
\end{proposition}

Notice that, by the remarks above, the converse of the proposition also
holds.

\begin{remark}
Equation (\ref{Fernandes:eqDirac}) can be interpreted in
terms of the Dirac theory of constraints. This is the reason for the
use of the term ``Poisson-Dirac submanifold''. We refer the reader to
\cite{Fernandes:CrFe} for more explanations.
\end{remark}

On the other hand, from Proposition \ref{Fernandes:propDirac}, we
deduce the following sufficient condition for a submanifold to be a
Poisson-Dirac submanifold:

\begin{corollary}
\label{Fernandes:corLieDirac}
Let $M$ be a Poisson manifold and $N\subset M$ a submanifold. Assume
that there exists a subbundle $E\subset T_NM$ such that: 
\[ T_N M=TN\oplus E\] 
and $\#(E^0)\subset TN$. Then $N$ is a Poisson-Dirac submanifold.
\end{corollary}

\begin{proof}
Under the assumptions of the corollary, one has a decomposition
\[ \pi=\pi_N+\pi_E,\]
where $\pi_N\in \Gamma(\wedge^2 TN)$ and $\pi_E\in \Gamma(\wedge^2
E)$ are both smooth bivector fields.  On the other hand, one checks
easily that (\ref{Fernandes:eqDirac}) holds. By Proposition
\ref{Fernandes:propDirac}, we conclude that $N$ is a Poisson-Dirac
submanifold.
\end{proof}

There are Poisson-Dirac submanifolds which do not satisfy the
conditions of this corollary. Also, the bundle $E$ may not be unique. 
For a detailed discussion and examples we refer to \cite{Fernandes:CrFe}. 

Under the assumptions of the corollary, the Poisson bracket on the
Poisson-Dirac submanifold $N\subset M$ is quite simple to describe:
Given two smooth functions $f,g\in C^\infty(N)$, to obtain their
Poisson bracket we pick extensions $\widetilde{f},\widetilde{g}\in
C^\infty(M)$ such that ${\mathrm d}_x\widetilde{f}, {\mathrm d}_x\widetilde{g}\in
E_x^0$. Then the Poisson bracket on $N$ is given by:
\begin{equation}
\label{Fernandes:eqPoissonBracket}
\left\{f,g\right\}_N={\{\widetilde{f},\widetilde{g}\}|}_N.
\end{equation}
It is not hard to check that this formula does not
depend on the choice of extensions.

\begin{remark}
\label{Fernandes:rmkLieDirac}
Let $M$ be a Poisson manifold and $N\subset M$ a submanifold. Assume
that there exists a subbundle $E\subset T_NM$ such that $E^0$ is a
\emph{Lie subalgebroid} of $T^*M$ (equivalently, $E$ is a co-isotropic
submanifold of the tangent Poisson manifold $TM$). Then $E$ satisfies
the assumptions of the corollary, so $N$ is a Poisson-Dirac submanifold.
This class of Poisson-Dirac submanifolds have very special geometric
properties. They where first study by Xu in \cite{Fernandes:Xu}, which
calls them \textbf{Dirac submanifolds}. They are further
discussed by Crainic and Fernandes in \cite{Fernandes:CrFe}, where
they are called \textbf{Lie-Dirac submanifolds}.
\end{remark}

\section{Fixed point sets of proper Poisson actions}

In this section we will give a proof of Theorem
\ref{Fernandes:MainTheorem}, which we restate now as follows:

\begin{theorem}
\label{Fernandes:MainTheoremAlt}
Let $G\times M\to M$ be a proper Poisson action. Then the fixed point
set $M^G$ is a Poisson-Dirac submanifold.
\end{theorem}

Since the action is proper, the fixed point set $M^G$ is an embedded
submanifold of $M$. Its connected components may have different
dimensions, but our argument will be valid for each such component, so
we will assume that $M^G$ is a connected submanifold. The proof will
consist in showing that there exists a subbundle $E\subset T_{M^G}M$
satisfying the conditions of Corollary \ref{Fernandes:corLieDirac}.

First of all, given any action $G\times M\to M$ (proper or not) there
exists a lifted action $G\times TM\to TM$. For proper actions we have
the following basic property:  

\begin{proposition}
If $G\times M\to M$ is a proper action then there exists a
$G$-invariant metric on $TM$.
\end{proposition}

For a proof of this fact and other elementary properties of proper
actions, we refer to \cite{Fernandes:DuKo}. Explicitly, the
$G$-invariance of the metric means that:
\[ 
\langle g\cdot v,g\cdot w\rangle_{g\cdot p}=\langle v,w\rangle_{p},
\qquad \forall v,w\in T_p M.
\]
where $g\in G$ and $p\in M$.

We fix, once and for all, a $G$-invariant metric $\langle~,~\rangle$
for our proper Poisson action $G\times M\to M$. Let us consider the
subbundle $E\subset T_{M^G}M$ which is orthogonal to $TM^G$:
\[ E=\{v\in T_{M^G}M:\langle v,w\rangle=0,\forall w\in TM^G\}.\]
We have:

\begin{lemma}
\[ T_{M^G} M=TM^G\oplus E\qquad\text{and}\qquad\#(E^0)\subset TM^G.\]
\end{lemma}

\begin{proof}
Since $E=(TM^G)^\perp$, the decomposition $T_{M^G}M=TM^G\oplus E$ is
obvious. Now for a proper action, we have $(TM)^G=TM^G$ so this
decomposition can also be written as:
\begin{equation}
\label{Fernandes:eqDecomp1}
T_{M^G} M=(TM)^G\oplus E,
\end{equation}
On the other hand, we have the lifted cotangent action $G\times
T^*M\to T^*M$, which is related to the lifted tangent action by
$g\cdot\xi(v)=\xi(g^{-1}\cdot v)$, $\xi\in T^*M, v\in TM$.
We claim that:
\begin{equation}
\label{Fernandes:eqDecomp2} 
E^0\subset (T^*M)^G.
\end{equation}
In fact, if $v\in TM$ we can use (\ref{Fernandes:eqDecomp1}) to decompose
it as $v=v_G+v_E$, where $v_G\in(TM)^G$ and $v_E\in E$. Hence, for
$\xi\in E^0$ we find: 
\begin{align*}
g\cdot\xi(v_G+v_E)
&=\xi(g^{-1}\cdot v_G+g^{-1}\cdot v_E)\\
&=\xi(v_G)+\xi(g^{-1}\cdot v_E)\\
&=\xi(v_G)\\
&=\xi(v_G)+\xi(v_E)=\xi(v_G+v_E).
\end{align*}
We conclude that $g\cdot\xi=\xi$ and (\ref{Fernandes:eqDecomp2}) follows.

Since $G\times M\to M$ is a Poisson action, we see that $\#:T^*M\to
TM$ is a $G$-equivariant bundle map. Hence, if $\xi\in E^0$, we obtain
from (\ref{Fernandes:eqDecomp2}) that:
\[ g\cdot\#\xi=\#(g\cdot\xi)=\#\xi.\]
This means that $\#\xi\in (TM)^G=TM^G$, so the lemma holds.
\end{proof}

This lemma shows that the conditions of Corollary \ref{Fernandes:corLieDirac} 
are satisfied, so $M^G$ is a Poisson-Dirac submanifold and the proof of
Theorem \ref{Fernandes:MainTheoremAlt} is completed.

\begin{remark}
{If} one works further with the decomposition
(\ref{Fernandes:eqDecomp1}) and its transposed version, it is not hard
to show that $E^0$ is actually a Lie subalgebroid of $T^*M$. Therefore,
the fixed point set $M^G$ of a proper Poisson action is, in fact, a
Lie-Dirac submanifold of $M$ (see Remark \ref{Fernandes:rmkLieDirac}).
\end{remark}

\begin{remark}
Special cases of Theorem \ref{Fernandes:MainTheoremAlt} where obtained
by Damianou and Fernandes in \cite{Fernandes:DaFe} for a compact Lie
group $G$, and by Fernandes and Vanhaecke in \cite{Fernandes:FeVan}
for a reductive algebraic group $G$. Xiang Tang also proves a version of
this theorem in his PhD thesis \cite{Fernandes:Tang}.
\end{remark}

Notice that the Poisson bracket of functions $f,g\in C^\infty(M^G)$
can be obtained simply by choosing $G$-invariant extensions
$\widetilde{f},\widetilde{g}\in C^\infty(M)^G$, and setting:
\[ \left\{f,g\right\}_{M^G}={\{\widetilde{f},\widetilde{g}\}|}_{M^G}.\]
This follows from equation (\ref{Fernandes:eqPoissonBracket}) and the
remark that for any such $G$-invariant extensions we have
${\mathrm d}_{M^G}\widetilde{f},{\mathrm d}_{M^G}\widetilde{g}\in E^0$. It is an instructive
exercise to prove directly that the bracket on ${M^G}$ does not depend on
the choice of extensions.

\section{Applications and further results}

Every compact Lie group action is proper. In particular,
a finite group action is always a proper. The case $G={\mathbb Z}_2$ leads to
the following result:

\begin{corollary}
Let $\phi:M\to M$ be an involutive Poisson automorphism of a Poisson
manifold $M$. The fixed point set $\{p\in M:\phi(p)=p\}$ has a natural
induced Poisson structure.
\end{corollary}

\begin{proof}
Apply Theorem \ref{Fernandes:MainTheoremAlt} to the Poisson action of 
the group $G=\{\text{Id},\phi\}$.
\end{proof}

This result is known in the literature as the Poisson Involution
Theorem (see \cite{Fernandes:DaFe,Fernandes:FeVan,Fernandes:Xu}). It
has been applied in \cite{Fernandes:DaFe,Fernandes:FeVan} to explain
the relationship between the geometry of the Toda and Volterra
lattices, and there should be similar relations between other known
integrable systems. In this respect, it should be interesting to find
extensions of our results to infinite dimensional manifolds and
actions.

\vskip 10 pt

Recall that if an action $G\times M\to M$ is proper and free then the
space of orbits $M/G$ is a smooth manifold. For general non-free
actions the orbit space can be a very pathological topological
space. However, for proper actions the singularities of the orbit
space are very much controlled, and $M/G$ is a nicely stratified
topological space. For proper symplectic actions there is a beautiful
theory of singular symplectic quotients due to Lerman and Sjamaar
\cite{Fernandes:LeSj} which describes the geometry of $M/G$. For
proper Poisson actions one should expect that the orbit space still
exhibits some nice Poisson geometry. In fact, we will explain in
\cite{Fernandes:FeOrRa} that Theorem \ref{Fernandes:MainTheoremAlt}
leads to the following result that generalizes a theorem due to Lerman
and Sjamaar:

\begin{theorem}
Let $G\times M\to M$ be a proper Poisson action. Then the quotient
$M/G$ is a Poisson stratified space.
\end{theorem}

Note that if a Poisson action is proper \emph{and} free then the orbit
space is a smooth Poisson manifold. In this case one can identify the
smooth functions on the quotient $M/G$ with the $G$-invariant
functions on $M$:
\[ C^\infty(M/G)\simeq C^\infty(M)^G.\]
In the non-free case, the smooth structure of $M/G$ as a stratified
space also leads to such an identification. Rather than explaining in
detail the notion of a Poisson stratified space (see the upcoming
paper \cite{Fernandes:FeOrRa}), we will illustrate this result with an
example.

\begin{example}
Let ${\mathbb C}^{n+1}$ be the complex $n+1$-dimensional space with holomorphic
coordinates $(z_0,\dots,z_n)$ and anti-holomorphic coordinates
$(\overline{z}_0,\dots,\overline{z}_n)$. On the (real) manifold 
${\mathbb C}^{n+1}-{0}$ we will consider a (real) quadratic Poisson bracket of the
form:
\[ \{z_i,z_j\}=a_{ij}z_iz_j,\qquad 
\{z_i,\overline{z}_j\}=\{\overline{z}_i,\overline{z}_j\}=0.\]
where $A=(a_{ij})$ is a skew-symmetric matrix.

The group ${\mathbb C}^*$ of non-zero complex numbers acts on ${\mathbb C}^{n+1}-0$ by
multiplication of complex numbers. This is a free and proper Poisson
action, so the quotient ${\mathbb CP}(n)={\mathbb C}^{n+1}-0/{\mathbb C}^*$ inherits a Poisson
bracket. 

Let us consider now the action of the $n$-torus ${\mathbb T}^n$ on ${\mathbb C}^{n+1}-0$
defined by:
\[ (\theta_1,\dots,\theta_n)\cdot (z_0,z_1,\cdots,z_n)=
(z_0,e^{i\theta_1}z_1,\cdots,e^{i\theta_n}z_n).\]
This is a Poisson action that commutes with the ${\mathbb C}^*$-action. It
follows that the ${\mathbb T}^n$-action descends to a Poisson action on
${\mathbb CP}(n)$. Note that the action of ${\mathbb T}^n$ on ${\mathbb CP}(n)$ is proper but
not free. The quotient ${\mathbb CP}(n)/{\mathbb T}^n$ is not a manifold but it can
be identified with the standard simplex 
\[ 
\Delta^n=\{(\mu_0,\dots,\mu_n)\in{\mathbb R}^{n+1}:\sum_{i=0}^n
\mu_i=1,\mu_i\ge 0\}.
\]
This identification is obtained via the map $\mu:{\mathbb CP}(n)\to\Delta^n$ defined by:
\[ 
\mu([z_0:\cdots:z_n])=
\left(\frac{|z_0|^2}{|z_0|^2+\cdots+|z_n|^2},\cdots,\frac{|z_n|^2}{|z_0|^2+\cdots+|z_n|^2}\right).
\]

Let us describe the Poisson stratification of
$\Delta^n={\mathbb CP}(n)/{\mathbb T}^n$. The Poisson bracket on $\Delta^n$ is obtained
through the identification:
\[ C^\infty(\Delta^n)\simeq C^\infty({\mathbb CP}(n))^{{\mathbb T}^n}.\]
For that, we simply compute the Poisson bracket between the components
of the map $\mu$. A more or less straightforward computation will
show that:
\begin{equation}
\label{Fernandes:eqPoissonBr}
\{\mu_i,\mu_j\}=\left(a_{ij}-\sum_{l=0}^n(a_{il}+a_{lj})\mu_l\right)\mu_i\mu_j,
\qquad (i,j=0,\dots,n).
\end{equation}
Now notice that (\ref{Fernandes:eqPoissonBr}) actually defines a Poisson
bracket on ${\mathbb R}^{n+1}$. For this Poisson bracket, the interior of the 
simplex and its faces are Poisson submanifolds: a face 
$\Delta_{i_1,\dots,i_{n-d}}$ of dimension $0\le d\le n$ is given by
equations of the form:
\[ \sum_{i=0}^n \mu_i=1,\quad \mu_{i_1}=\cdots=\mu_{i_{n-d}}=0,
\quad\mu_i>0 \text{ for }i\not\in\{i_1,\dots,i_{n-d}\}.
\]
These equations define Poisson submanifolds since:
\begin{enumerate}
\item[(a)] the bracket $\{\mu_i,\mu_l\}$ vanishes whenever $\mu_l=0$, and
\item[(b)] the bracket $\{\mu_i,\sum_{l=0}^n\mu_l\}$ vanishes whenever 
$\sum_{l=0}^n\mu_l=1$.
\end{enumerate}
Therefore, the Poisson stratification of $\Delta^n$ consists of
strata formed by the faces of dimension $0\le d\le n$, which are
smooth Poisson manifolds.
\end{example}


\begin{thebibliography}{99}

\bibitem{Fernandes:CrFe} 
M.~Crainic and R.L.~Fernandes, Integrability of Poisson brackets,
\emph{Journal of Differential Geometry} \textbf{66} (2004), 71--137.

\bibitem{Fernandes:DaFe}
P.~Damianou and R.L.~Fernandes, From the Toda lattice to the Volterra
lattice and back, \emph{Reports on Math.~Phys.~}\textbf{50},(2002)
361--378.

\bibitem{Fernandes:DuKo} J.~Duistermaat and J.~Kolk, \emph{Lie Groups},
  Springer-Verlag Berlin Heidelberg, 2000.

\bibitem{Fernandes:FeVan}
R.L.~Fernandes and P.~Vanhaecke, Hyperelliptic Prym Varieties and
Integrable Systems, \emph{Commun.~Math.~Phys.}~\textbf{221} (2001)
169--196.

\bibitem{Fernandes:FeOrRa}
R.L.~Fernandes, J.-P.~Ortega and T.~Ratiu, Momentum maps in Poisson
geometry, paper in preparation.

\bibitem{Fernandes:GuSt}
V.~Guillemin and S.~Sternberg, Convexity properties of the moment
mapping, \emph{Inventiones Mathematicae} \textbf{67} (1982), 491--513.

\bibitem{Fernandes:OrRa}
J.-P.~Ortega and T.~Ratiu, \emph{Momentum maps and {H}amiltonian reduction},
Progress in Mathematics, vol.~222, Birkh\"auser, Boston, 2004.

\bibitem{Fernandes:LeSj}
Eugene Lerman and Reyer Sjamaar, Stratified symplectic spaces
and reduction, \emph{Annals of Mathematics (2)} \textbf{134} (1991),
375--422.

\bibitem{Fernandes:Tang}
X.~Tang, \emph{Quantization of Noncommutative Poisson manifolds}, PhD Thesis,
University of California, Berkeley, USA, (2004).

\bibitem{FernandesVai}
I.~Vaisman, \emph{Lectures on the Geometry of Poisson Manifolds},
Progress in Mathematics, vol.~\textbf{118}, Birkh\"auser, Berlin, 1994.

\bibitem{Fernandes:Xu}
P.~Xu, Dirac submanifolds and Poisson involutions, 
\emph{Ann.~Sci.~\'Ecole Norm.~Sup.~}(4) \textbf{36} (2003), 403--430.

\end{thebibliography}
\end{document}